\input amstex
\input amsppt.sty
\magnification=\magstep1
\hsize=33truecc
\vsize=22.2truecm
\baselineskip=16truept
\NoBlackBoxes
\TagsOnRight \pageno=1 \nologo
\def\Z{\Bbb Z}
\def\N{\Bbb N}

\def\l{\left}
\def\r{\right}
\def\bg{\bigg}
\def\({\bg(}
\def\[{\bg\lfloor}
\def\){\bg)}
\def\]{\bg\rfloor}
\def\t{\text}
\def\f{\frac}

\def\bi{\binom}
\def\eq{\equiv}

\def\ls{\leqslant}

\def\mo{\roman{mod}}

\def\al{\alpha}

\def\FF#1#2#3{{}_2F_1\bigg(\matrix\format\c\\#1\\#2\endmatrix\bigg|#3\bigg)}
\def\FFF#1#2#3{{}_3F_2\bigg(\matrix\format\c\\#1\\#2\endmatrix\bigg|#3\bigg)}

\def\Proof{\noindent{\it Proof}}

\def\Remark{\medskip\noindent{\it  Remark}}

\hbox {Preprint, {\tt arXiv:1210:5237}}
\bigskip
\topmatter
\title Determining $x$ or $y$ mod $p^2$ with $p=x^2+dy^2$ \endtitle
\author Zhi-Wei Sun\endauthor
\leftheadtext{Zhi-Wei Sun}
\rightheadtext{Determining $x$ or $y$ mod $p^2$ with $p=x^2+dy^2$}
\affil Department of Mathematics, Nanjing University\\
 Nanjing 210093, People's Republic of China
  \\  zwsun\@nju.edu.cn
  \\ {\tt http://math.nju.edu.cn/$\sim$zwsun}
\endaffil
\abstract Let $p$ be an odd prime and let $d\in\{2,3,7\}$. When $(\f{-d}p)=1$ we can write $p=x^2+dy^2$
with $x,y\in\Z$; in this paper we aim at determining $x$ or $y$ modulo $p^2$. For example,
when $p=x^2+3y^2$, we show that if $p\eq x\eq 1\ (\mo\ 4)$ then
$$\sum_{k=0}^{(p-1)/2}(3[3\mid k]-1)(2k+1)\f{\bi{2k}k^2}{(-16)^k}\eq\l(\f2p\r)2x\pmod{p^2}$$
where $[3\mid k]$ takes $1$ or $0$ according as $3\mid k$ or not, and that if $-p\eq y\eq 1\ (\mo\ 4)$ then
$$\sum_{k=0}^{(p-1)/2}\l(\f k3\r)\f{k\bi{2k}k^2}{(-16)^k}
\eq(-1)^{(p+1)/4}y\eq\sum_{k=0}^{(p-1)/2}(1-3[3\mid k])\f{k\bi{2k}k^2}{(-16)^k}\pmod{p^2}.$$
We also determine 
$$\sum_{k=0}^{p-1}\f{k\bi{2k}k^3}{m^k}\sum_{k\ls j<2k}\f1j\quad \mo\ p$$
for $m=1,-8,16,-64,256,-512,4096$.
\endabstract
\thanks 2010 {\it Mathematics Subject Classification}.\,Primary 11E25, 11A07;
Secondary  05A10, 11B65.
\newline\indent Supported by the National Natural Science
Foundation (Grant No. 11171140) of China.
\endthanks
\endtopmatter
\document

\heading{1. Introduction}\endheading

Let $p\eq1\pmod4$ be a prime. It is well-known that we may write $p=x^2+y^2$ with
$x\eq1\pmod4$ and $y\eq0\pmod2$. In 1828 Gauss determined $x$ mod
$p$ via the congruence
$$\bi{(p-1)/2}{(p-1)/4}\eq2x\pmod p$$
(see, e.g., [BEW, Chapter 9]).
In 1986, S. Chowla, B. Dwork and R. J. Evans [CDE] showed further that
$$\bi{(p-1)/2}{(p-1)/4}\eq\f{2^{p-1}+1}2\l(2x-\f p{2x}\r)\pmod {p^2},$$
which implies the congruence
$$\bi{(p-1)/2}{(p-1)/4}^2\eq 2^{p-1}(4x^2-2p)\pmod{p^2}.$$

Let $p$ be an odd prime. In 1997 van Hamme [vH] showed that
$$\sum_{k=0}^{p-1}\f{\bi{2k}k^3}{64^k}\eq\cases4x^2-2p\ (\mo\ p^2)&\t{if}\ p\eq1\ (\mo\ 4)\ \&\ p=x^2+y^2\ (2\nmid x),
\\0\ (\mo\ p^2)&\t{if}\ p\eq3\ (\mo\ 4).\endcases\tag1.0$$
By Stirling's formula, $n!\sim\sqrt{2\pi n}(n/e)^n$ and hence
$$\lim_{k\to\infty}\root k\of{\bi{2k}k^3}=64.$$
Thus we call $\sum_{k=0}^{p-1}\bi{2k}k^3/m^k$ mod $p^2$ the {\it
border case} if $m=64$. In 2009 the author posed several
conjectures on $\sum_{k=0}^{p-1}\bi{2k}k^3/m^k$ mod $p^2$ in
non-border cases. For example, in Nov. 2009 he conjectured that
$$\aligned\sum_{k=0}^{p-1}\bi{2k}k^3\eq\cases4x^2-2p\ (\mo\ p^2)&\t{if}\ (\f p7)=1\ \&\ p=x^2+7y^2,
\\0\ (\mo\ p^2)&\t{if}\ (\f p7)=-1.\endcases\endaligned\tag1.1$$
This appeared as Conjecture 5.3 of [Su1]. In [Su2, Conjecture 5.2]
the author also formulated conjectures on
$\sum_{k=0}^{p-1}\bi{2k}k^3/m^k$ with $m=-8,16,-64,256,-512,4096$.
Namely,
$$\sum_{k=0}^{p-1}\f{\bi{2k}k^3}{(-8)^k}\eq\sum_{k=0}^{p-1}\f{\bi{2k}k^3}{(-512)^k}
\eq\sum_{k=0}^{p-1}\f{\bi{2k}k^3}{64^k}\ \l(\mo\ p^{(5+(\f{-1}p)/2)}\r),\tag1.2$$
$$\sum_{k=0}^{p-1}\f{\bi{2k}k^3}{(-64)^k}\eq\cases4x^2-2p\ (\mo\ p^2)&\t{if}\  p=x^2+2y^2,
\\0\ (\mo\ p^2)&\t{if}\ p\eq5,7\ (\mo\ 8),\endcases\tag1.3$$
$$\sum_{k=0}^{p-1}\f{\bi{2k}k^3}{16^k}\eq\cases4x^2-2p\ (\mo\ p^2)&\t{if}\  p=x^2+3y^2,
\\0\ (\mo\ p^2)&\t{if}\ p\eq2\ (\mo\ 3),\endcases\tag1.4$$
$$\sum_{k=0}^{p-1}\f{\bi{2k}k^3}{256^k}\eq\l(\f{-1}p\r)\sum_{k=0}^{p-1}\f{\bi{2k}k^3}{16^k}\ \l(\mo\ p^{(5+(\f p3))/2}\r),\tag1.5$$
and
$$\sum_{k=0}^{p-1}\f{\bi{2k}k^3}{4096^k}\eq\l(\f{-1}p\r)\sum_{k=0}^{p-1}\bi{2k}k^3\ \l(\mo\ p^{(5+(\f p7))/2}\r).\tag1.6$$
Note that [A, Theorem 5] implies (1.1)-(1.6) modulo $p$.
In an earlier preprint form of [S2] posted to arXiv, Z.-H. Sun got partial success in proving the above conjectures; in fact, he showed that
$$\gather\sum_{k=0}^{p-1}\bi{2k}k^3\eq\sum_{k=0}^{p-1}\f{\bi{2k}k^3}{4096^k}\eq0\ (\mo\ p^2)\ \ \t{if}\ \l(\f p7\r)=-1,
\\\sum_{k=0}^{p-1}\f{\bi{2k}k^3}{(-8)^k}\eq\sum_{k=0}^{p-1}\f{\bi{2k}k^3}{(-512)^k}\eq 0\ (\mo\ p^2)\ \ \t{if}\ p\eq3\ (\mo\ 4),
\\\sum_{k=0}^{p-1}\f{\bi{2k}k^3}{16^k}\eq\sum_{k=0}^{p-1}\f{\bi{2k}k^3}{256^k}\eq 0\ (\mo\ p^2)\ \ \t{if}\ p\eq2\ (\mo\ 3),
\endgather$$
and
$$\sum_{k=0}^{p-1}\f{\bi{2k}k^3}{(-8)^k}\eq 4x^2-2p\ (\mo\ p^2)\quad\t{if}\ p=x^2+y^2\ (2\nmid x).$$
The key tool of Z.-H. Sun [S2] is the following identity
$$P_n(x)^2=\sum_{k=0}^n\bi nk\bi{n+k}k\bi{2k}k\l(\f{x^2-1}4\r)^k,\tag1.7$$
where $P_n(x)$ is the Legendre polynomial of degree $n$ given by
$$P_n(x)=\sum_{k=0}^n\bi nk\bi{n+k}k\l(\f{x-1}2\r)^k.$$
Note that (1.7) follows from the well-known Clausen identity
$$\FF{2a,2b}{a+b+1/2}{z}^2=\FFF{2a,2b,a+b}{a+b+1/2,2a+2b}{4z(1-z)}$$
in the special case $a=-n/2$, $b=(n+1)/2$ and $z=(1-x)/2$.
In 2012, the author's conjectural congruences (1.1)-(1.6)
modulo $p^2$ were finally confirmed by J. Kibelbek, L. Long, K. Moss, B.
Sheller and H. Yuan in the preprint [KLMSY] posted to arXiv in Oct. 2012,
as well as the final version of [S2].
Both [KLMSY] and [S2] depend heavily on the paper [CV] of M. J. Coster and L. van Hamme 
in which the authors determined
$$P_n(\sqrt2),\ P_n\l(\f{3\sqrt2}4\r),\ P_n(\sqrt{-3}),\ P_n\l(\f{\sqrt3}2\r),\ P_n\l(3\sqrt{-7}\r),\ P_n\l(\f{3\sqrt7}8\r)$$
modulo $p^2$ via complex multiplication of elliptic curves, where $p=2n+1$ is an odd prime.

 Concerning (1.4) we mention that the author [Su3] proved that
 $$\sum_{k=0}^{p-1}(-1)^kA_k\eq\sum_{k=0}^{p-1}\f{\bi{2k}k^3}{16^k}\pmod{p^2}\tag1.8$$
 for any odd prime $p$, where $A_0,A_1,A_2,\ldots$ are Ap\'ery numbers given by
 $$A_n=\sum_{k=0}^n\bi nk^2\bi{n+k}k^2\ (n=0,1,2,\ldots).$$
 The reader may consult [Su7] and [Su8] for many other conjectures similar to (1.1)-(1.6).

 During their initial attempt to prove (1.1) in 2010, M. Jameson and K. Ono [JO]
realized that $$\sum_{k=0}^{(p-1)/2}\bi{2k}k^3(H_{2k}-H_k)\eq0\pmod{p}\quad \t{for any prime}\ p>3$$
but they did not have a proof of this observation, where $H_k$ denotes the harmonic number $\sum_{0<j\ls k}1/j$.
When $p>3$ is a prime with $p\eq 3\ (\mo\ 4)$, by [Su1, (1.11)] and (1.0) we have
$$\sum_{k=0}^{(p-1)/2}\f{\bi{2k}k^3}{64^k}H_{k}\eq0\ (\mo\ p).$$
On March 2, 2010 the author formulated a conjecture (which was
presented as [Su, Conj. A37]) on
$\sum_{k=0}^{p-1}\bi{2k}k^3(H_{2k}-H_k)/m^k$ modulo $p$ or $p^2$,
where $p>3$ is a prime and $m\in\{1,-8,16,-64,256,-512,4096\}$.
In [KLMSY] it was proved for any prime $p>3$
and $m=1,-8,16,-64,256,-512,4096$ the congruence
$$\sum_{k=0}^{p-1}\f{\bi{2k}k^3}{m^k}(H_{2k}-H_k)\eq \f{q_p(m)}6\sum_{k=0}^{p-1}\f{\bi{2k}k^3}{m^k}\ (\mo\ p)\tag1.9$$
where $q_p(m)$ denotes the Fermat quotient $(m^{p-1}-1)/p$.

When $p=x^2+y^2$ with $4\mid x-1$ and $2\mid y$, the author [Su5] showed that we can determine $x$ mod $p^2$ in the following way:
$$\l(\f 2p\r)x\eq\sum_{k=0}^{(p-1)/2}\f{k+1}{8^k}\bi{2k}k^2\eq\sum_{k=0}^{(p-1)/2}\f{2k+1}{(-16)^k}\bi{2k}k^2\pmod{p^2}.$$

With the help of Coster-van Hamme's work, in this paper we determine $x$ or $y$ mod $p^2$ when $p=x^2+dy^2$ with $d=2,3,7$.

Before stating our main results we need to introduce Lucas sequences.

Let $A,B\in\Z$. The Lucas sequence $u_n=u_n(A,B)\ (n\in\N)$ and its companion $v_n=v_n(A,B)\ (n\in\N)$  are defined as follows:
$$\gather u_0=0,\ u_1=1,\ \t{and}\ u_{n+1}=Au_n-Bu_{n-1}\ (n=1,2,3,\ldots);
\\ v_0=2, v_1=A,\  \t{and}\ v_{n+1}=Av_n-Bv_{n-1}\ (n=1,2,3,\ldots).\endgather$$
It is well known that $(\al-\beta)u_n=\al^n-\beta^n$ and $v_n=\al^n+\beta^n$ for all $n\in\N$, where $\al$ and $\beta$ are the two roots
of the equation $x^2-Ax+B=0$. The sequence $P_n=u_n(2,-1)\ (n\in\N)$ is called the Pell sequence and its companion is the sequence
$Q_n=v_n(2,1)\ (n\in\N)$.

\proclaim{Theorem 1.1}  Let $p$ be a prime
with $p\eq1,3\ (\mo\ 8)$.
And write $p=x^2+2y^2$ with $x,y\in\Z$ and $x\eq1\ (\mo\ 4)$.

{\rm (i)} If $p\eq1\ (\mo\ 8)$, then
$$\align&4\sum_{k=0}^{p-1}\bi{2k}k^2\f{kP_k}{32^k}\eq\sum_{k=0}^{p-1}\bi{2k}k^2\f{kQ_k}{32^k}
\\\eq&(-1)^{(p-1)/8+(x-1)/4}\l(\f px-2x\r)\ (\mo\ p^2),\endalign$$
and we can determine $x$ mod $p^2$ via the congruence
$$(-1)^{(x-1)/4}x\eq\f{(-1)^{(p-1)/8}}2\sum_{k=0}^{p-1}\bi{2k}k^2\f{(k+1)Q_k}{32^k}\ (\mo\ p^2).$$

{\rm (ii)} If $p\eq3\ (\mo\ 8)$, then
$$\sum_{k=0}^{p-1}\bi{2k}k^2\f{kP_k}{32^k}\eq\f12\sum_{k=0}^{p-1}\bi{2k}k^2\f{kQ_k}{32^k}\eq(-1)^{(y+1)/2}y\ (\mo\ p^2)$$
and
$$\sum_{k=0}^{p-1}\bi{2k}k^2\f{P_k}{32^k}\eq (-1)^{(y-1)/2}\l( 2y-\f p{4y}\r)\ (\mo\ p^2).$$
\endproclaim
\Remark\ 1.1. Theorem 1.1 was first conjectured by the author in 2009, it appeared as part of [Su4, Conjecture 2.3].

\proclaim{Theorem 1.2} Let $p\eq 1\ (\mo\ 3)$ be a prime and write $p=x^2+3y^2$ with $x,y\in\Z$.

{\rm (i)} If $p\eq x\eq1\ (\mo\ 4)$, then we can determine $x$ mod $p^2$ via the congruence
$$\l(\f 2p\r) 2x\eq\sum_{k=0}^{p-1}(3[3\mid k]-1)(2k+1)\f{\bi{2k}k^2}{(-16)^k}\pmod{p^2}.\tag1.10$$

{\rm (ii)} If $-p\eq y\eq1\ (\mo\ 4)$, then
$$\sum_{k=0}^{p-1}\l(\f k3\r)\f{\bi{2k}k^2}{(-16)^k}\eq(-1)^{(p-3)/4}\l(4y-\f p{3y}\r)\pmod{p^2},\tag1.11$$
and we can determine $y$ mod $p^2$ in the following way:
$$\sum_{k=0}^{p-1}\l(\f k3\r)\f{k\bi{2k}k^3}{(-16)^k}\eq(-1)^{(p+1)/4}y\eq\sum_{k=0}^{p-1}(1-3[3\mid k])\f{k\bi{2k}k^2}{(-16)^k}\pmod{p^2}.\tag1.12$$
\endproclaim
\Remark\ 1.2. Part (ii) of Theorem 1.2 appeared as part (i) of [Su2, Conjecture 5.11].

\proclaim{Theorem 1.3} Let $p$ be an odd prime with $(\f p7)=1$ and write $p=x^2+7y^2$ with $x,y\in\Z$.

{\rm (i)} If $p\eq x\eq1\ (\mo\ 4)$, then we may determine $x$ mod $p^2$ via the congruence
$$\sum_{k=0}^{p-1}(4k+3)\f{\bi{2k}k^2}{16^k}v_k(1,16)\eq 6\l(\f 2p\r) x\pmod{p^2}.\tag1.13$$

{\rm (ii)} If $-p\eq y\eq1\ (\mo\ 4)$, then we may determine $y$ mod $p^2$ in the following way:
$$\sum_{k=0}^{p-1}\f{k\bi{2k}k^2}{16^k}u_k(1,16)\eq\sum_{k=0}^{p-1}\f{k\bi{2k}k^2}{16^k}v_k(1,16)\eq-\l(\f 2p\r)\f y2\pmod{p^2}.\tag1.14$$
\endproclaim

\proclaim{Theorem 1.4} Let $p>3$ be a prime.

{\rm (i)} If $p\eq1\ (\mo\ 12)$ and $p=x^2+3y^2$ with $x,y\in\Z$ and $x\eq1\ (\mo\ 4)$, then
$$\sum_{k=0}^{p-1}\bi{2k}k^2\f{v_k(4,1)}{64^k}\eq4x-\f p{x}\ (\mo\ p^2)$$
and we can determine $x$ mod $p^2$ by
$$\sum_{k=0}^{p-1}\bi{2k}k^2\f{(k-1)v_k(4,1)}{64^k}\eq -2x\ (\mo\ p^2).$$

{\rm (ii)} If $p\eq7\ (\mo\ 12)$ and $p=x^2+3y^2$ with $x,y\in\Z$ and $y\eq1\ (\mo\ 4)$, then
$$\sum_{k=0}^{p-1}\f{u_k(4,1)}{64^k}\bi{2k}k^2\eq 2y-\f p{6y}\ (\mo\ p^2)$$
and also
$$y\eq\sum_{k=0}^{p-1}\bi{2k}k^2\f{ku_k(4,1)}{64^k}\eq\f14\sum_{k=0}^{p-1}\bi{2k}k^2\f{kv_k(4,1)}{64^k}
\ (\mo\ p^2).$$
\endproclaim
\Remark\ 1.3. Theorem 1.4 was part of [Su4, Conjecture 2.4].
\medskip

In contrast with (1.9), we have the following theorem.

\proclaim{Theorem 1.5} Let $p$ be an odd prime and let $a,b,m\in\Z$ with $m\not\eq0\ (\mo\ p)$.
Then
$$\aligned&\l(\f{-m}p\r)\sum_{k=0}^{p-1}\f{ak+b}{m^k}\bi{2k}k^3+\f{\bar m^{p-1}+1}4\sum_{k=0}^{p-1}\f{2ak+a-2b}{\bar m^k}\bi{2k}k^3
\\\eq&\f {ap}2\sum_{k=0}^{p-1}\f{\bi{2k}k^3}{\bar m^k}+\f{3p}2\sum_{k=0}^{p-1}\f{2ak+a-2b}{\bar m^k}\bi{2k}k^3(H_{2k}-H_k)\pmod{p^2},
\endaligned\tag1.15$$
where $\bar m=2^{12}/m$.
\endproclaim

Clearly Theorem 1.5 with $a=0$ and $b=1$ yields the following corollary.

\proclaim{Corollary 1.1} Let $p$ be an odd prime and let $m\in\Z$ with $p\nmid m$.
Set $\bar m=4096/m$. Then
$$\aligned&\l(\f{-m}p\r)\sum_{k=0}^{p-1}\f{\bi{2k}k^3}{m^k}-\sum_{k=0}^{p-1}\f{\bi{2k}k^3}{\bar m^k}
\\\eq&-\f p2\sum_{k=0}^{p-1}\f{\bi{2k}k^3}{\bar m^k}(6(H_{2k}-H_k)-q_p(\bar m))\pmod{p^2}.
\endaligned$$
\endproclaim

\Remark\ 1.4. In view of Corollary 1.1, for any odd prime $p>3$, (1.9) holds for all $m\in\{1,-8,16,-64,256,-512,4096\}$,
if and only if the symmetric relation
$$\l(\f{-m}p\r)\sum_{k=0}^{p-1}\f{\bi{2k}k^3}{m^k}\eq \sum_{k=0}^{p-1}\f{\bi{2k}k^3}{\bar m^k}\pmod{p^2}\tag1.16$$
holds for all $m\in\{1,-8,16,-64,256,-512,4096\}$. The latter was conjectured by the author in [Su2, Conj. 5.2].

In the next section shall provide some lemmas. We are going to show Theorems 1.1-1.4 in Section 3.
In Section 4 we will prove Theorem 1.5 and derive several corollaries of Th. 1.5. Section 4 also contains a conjecture
involving harmonic numbers.

\heading{2. Several lemmas}\endheading

\proclaim{Lemma 2.1} Let $p$ be an odd prime and let $m\in\Z$ with $m\not\eq0\ (\mo\ p)$.
Let $m_*$ be a root of the equation $x^2-mx+16m=0$ in the ring of algebraic $p$-adic integers.
Then, for any algebraic $p$-adic integers $a$ and $b$, we have
$$\aligned&\sum_{k=0}^{p-1}\f{\bi{2k}k^3}{m^k}\l(\f{ak}{16}(m_*-m+32)+b\r)
\\\eq&2a\sum_{k=0}^{p-1}\f{k\bi{2k}k^2}{m_*^k}\sum_{j=0}^{p-1}\f{\bi{2j}j^2}{m_*^j}+b\(\sum_{k=0}^{p-1}\f{\bi{2k}k^2}{m_*^k}\)^2\pmod{p^2}.
\endaligned\tag2.1$$
\endproclaim
\Proof. Set
$$S:=\sum_{n=0}^{p-1}\f{an+b}{m_*^n}\sum_{k=0}^n\bi{2k}k^2\bi{2(n-k)}{n-k}^2.$$
By [Su5, Lemma 3.1],
$$\align S=&\sum_{n=0}^{p-1}\f{an+b}{m_*^n}\sum_{k=0}^n\bi{2k}k^3\bi{k}{n-k}(-16)^{n-k}
\\=&\sum_{k=0}^{p-1}\f{\bi{2k}k^3}{m_*^k}\sum_{n=k}^{p-1}(an+b)\bi{k}{n-k}\l(-\f{16}{m_*}\r)^{n-k}
\\\eq&\sum_{k=0}^{(p-1)/2}\f{\bi{2k}k^3}{m_*^k}\sum_{j=0}^k(a(k+j)+b)\bi kj\l(-\f{16}{m_*}\r)^j\pmod{p^2}.
\endalign$$
Thus
$$\align S\eq&\sum_{k=0}^{(p-1)/2}\f{\bi{2k}k^3}{m_*^k}\((ak+b)\l(1-\f{16}{m_*}\r)^k+ak\sum_{0<j\ls k}\bi{k-1}{j-1}\l(-\f{16}{m_*}\r)^j\)
\\=&\sum_{k=0}^{(p-1)/2}\f{\bi{2k}k^3}{m_*^k}\((ak+b)\l(1-\f{16}{m_*}\r)^k-ak\f{16}{m_*}\l(1-\f{16}{m_*}\r)^{k-1}\)
\\=&\sum_{k=0}^{(p-1)/2}\f{\bi{2k}k^3}{m_*^k}\l(1-\f{16}{m_*}\r)^k\l(ak+b-\f{16}{m_*}ak\f m{m_*}\r)
\\=&\sum_{k=0}^{(p-1)/2}\f{\bi{2k}k^3}{m^k}\l(\f{ak}{16}(m_*-m+32)+b\r)\pmod{p^2}.
\endalign$$
On the other hand,
$$\align S=&\sum_{k=0}^{p-1}\f{\bi{2k}k^2}{m_*^k}\sum_{n=k}^{p-1}(an+b)\f{\bi{2(n-k)}{n-k}^2}{m_*^{n-k}}
\\\eq&\sum_{k=0}^{(p-1)/2}\f{\bi{2k}k^2}{m_*^k}\sum_{j=0}^{(p-1)/2}(a(k+j)+b)\f{\bi{2j}j^2}{m_*^j}
\\=&2a\sum_{k=0}^{p-1}\f{k\bi{2k}k^2}{m_*^k}\sum_{j=0}^{p-1}\f{\bi{2j}j^2}{m_*^j}+b\(\sum_{k=0}^{p-1}\f{\bi{2k}k^2}{m_*^k}\)^2\pmod{p^2}.
\endalign$$
Combining the above, we obtain (2.1). \qed

\proclaim{Lemma 2.2} Let $p=2n+1$ be an odd prime. Then
$$P_n(\pm x)\eq(\pm1)^n\sum_{k=0}^{p-1}\f{\bi{2k}k^2}{(-16)^k}\l(\f{x-1}2\r)^k\pmod{p^2}.\tag2.2$$
\endproclaim
\Proof. It is well known that $P_n(-x)=(-1)^nP_n(x)$. As noted by van Hammer [vH], for each $k=0,\ldots,n$ we have
$$\align \bi nk\bi{n+k}k=&\bi nk\bi{-n-1}k(-1)^k=\bi{(p-1)/2}k\bi{(-p-1)/2}k(-1)^k
\\\eq&\bi{-1/2}k^2(-1)^k=\f{\bi{2k}k^2}{(-16)^k}\pmod{p^2}.
\endalign$$
So (2.2) holds. \qed

\proclaim{Lemma 2.3} Let $d\in\{1,2,3,7\}$ and suppose that $p$ is an odd prime with $(\f{-d}p)=-1$.
Let $\sqrt{-d}$ be a solution of $z^2+d=0$ in the ring of $p$-adic integers. Write $p=x^2+dy^2$ with $x,y\in\Z$
such that $\pi=x+y\sqrt{-d}\eq0\ (\mo\ p)$. Then, for $\bar\pi=x-y\sqrt{-d}$ we have
$$\bar\pi\eq 2x-\f p{2x}\eq -\f{\sqrt{-d}}2\l(4y-\f p{dy}\r)\pmod{p^2}.\tag2.3$$
\endproclaim
\Proof.  Note that
$$\bar\pi =2x-(x+y\sqrt{-d})=2x-\f p{x-\sqrt{-d}y}\eq 2x-\f p{2x}\pmod{p^2}$$
and
$$\bar\pi=-2y\sqrt{-d}+\f p{x-y\sqrt{-d}}\eq-2y\sqrt{-d}+\f p{-2y\sqrt{-d}}=-\f{\sqrt{-d}}2\l(4y-\f p{dy}\r)\ (\mo\ p^2).$$
This concludes the proof. \qed

\proclaim{Lemma 2.4} Let $p=2n+1$ be an odd prime. Let $d\in\{2,3,7\}$ and suppose that $(\f{-d}p)=1$.
Write $p=x^2+dy^2$ such that $\pi=x+y\sqrt{-d}\eq0\ (\mo\ p)$, and that  $x\eq1\ (\mo\ 4)$ if $d=2$, and $x+y\eq1\ (\mo\ 4)$ if $d\in\{3,7\}$.
Set $\bar\pi=x-y\sqrt{-d}$.

{\rm (i)} If $d=2$, then
$$P_n(\sqrt2)\eq i^{-y}\bar \pi\pmod{p^2}.$$

{\rm (ii)} If $d=3$, then
$$P_n(\sqrt{-3})\eq(-1)^y\bar\pi\pmod{p^2}$$
and
$$P_n\l(\f{\sqrt 3}2\r)\eq(-i)^n\bar\pi\pmod{p^2}.$$

{\rm (iii)} If $d=7$, then
$$P_n(3\sqrt{-7})\eq(-1)^{n+y}\bar\pi\pmod{p^2}$$
and
$$P_n\l(\f{3\sqrt7}8\r)\eq i^n\bar\pi\pmod{p^2}.$$
\endproclaim
\Proof. This follows from [CV, (48)] in the case $m=r=1$ and Table III of [CV, p.\,282].
However, we note that certain entries of Table III should be corrected. This results in our addition of $(-1)^y$
in part (ii) and the first congruence in part (iii). (One can check that our correction is necessary via computation.) \qed

\heading{3. Proofs of Theorems 1.1-1.4}\endheading

\medskip
\noindent{\it Proof of Theorem 1.1}. Without loss of generality we assume that $\pi=x+y\sqrt{-2}\eq0\ (\mo\ p)$.
(If not, we may use $-y$ instead of $y$.) It is easy to see that $(-1)^y=(-1)^n$. By Lemmas 2.2 and 2.4,
$$\sum_{k=0}^{p-1}\f{\bi{2k}k^2}{32^k}(1\pm\sqrt2)^k\eq (\mp1)^nP_n(\sqrt2)\eq (\mp1)^n(-i)^y\bar\pi=(\pm1)^ni^y\bar\pi\pmod{p^2}.\tag3.1$$
Thus
$$\align 2\sqrt2\sum_{k=0}^{p-1}\f{\bi{2k}k^2}{32^k}P_k=&\sum_{k=0}^{p-1}\f{\bi{2k}k^2}{32^k}((1+\sqrt2)^k-(1-\sqrt2)^k)
\\\eq&i^y\bar\pi(1-(-1)^n)\eq\cases0\ (\mo\ p^2)&\t{if}\ p\eq1\ (\mo\ 8),
\\2i^y\bar\pi\ (\mo\ p^2)&\t{if}\ p\eq 3\ (\mo\ 8).\endcases
\endalign$$
By Lemma 2.3,
$$\bar\pi\eq-\sqrt{-2}\l(2y-\f p{4y}\r)\pmod{p^2}.$$
So, if $p\eq3\ (\mo\ 8)$ then
$$2\sum_{k=0}^{p-1}\f{\bi{2k}k^2}{32^k}P_k\eq i^y\sqrt2\bar\pi\eq(-2i)i^y\l(2y-\f p{4y}\r)=2i^{y-1}\l(2y-\f p{4y}\r)\pmod{p^2}$$
and hence
$$\sum_{k=0}^{p-1}\f{\bi{2k}k^2}{32^k}P_k\eq(-1)^{(y-1)/2}\l(2y-\f p{4y}\r)\pmod{p^2}$$
as desired.
Note that
$$\align \sum_{k=0}^{p-1}\f{\bi{2k}k^2}{32^k}Q_k=&\sum_{k=0}^{p-1}\f{\bi{2k}k^2}{32^k}((1+\sqrt2)^k+(1-\sqrt2)^k)
\\\eq&i^y\bar\pi(1+(-1)^n)\eq\cases2i^y\bar\pi\ (\mo\ p^2)&\t{if}\ p\eq1\ (\mo\ 8),
\\0\ (\mo\ p^2)&\t{if}\ p\eq 3\ (\mo\ 8).\endcases
\endalign$$
If $p\eq1\ (\mo\ 8)$, then $2\mid y$ and
$$\f y2\eq\l(\f y2\r)^2+\f{x-1}4\cdot\f{x-1}2=\f{y^2}4+\f{x^2-1}8-\f{x-1}4=\f{p-1}8-\f{x-1}4\pmod2,$$
hence
$$\sum_{k=0}^{p-1}\f{\bi{2k}k^2}{32^k}Q_k\eq2(-1)^{y/2}\bar\pi\eq(-1)^{(p-1)/8+(x-1)/4}\l(4x-\f px\r)\pmod{p^2}$$
with the help of Lemma 2.3.

Applying Lemma 2.1 with $m=-64$, $m_*=-32\pm32\sqrt2$, $a=2$ and $b=2\pm\sqrt2$, we obtain that
$$\align&(2\pm\sqrt2)\sum_{k=0}^{p-1}(4k+1)\f{\bi{2k}k^3}{(-64)^k}
\\\eq&
4\sum_{k=0}^{p-1}\f{k\bi{2k}k^3}{32^k}(1\pm\sqrt2)^k\sum_{j=0}^{p-1}\f{\bi{2j}j^2}{32^j}(1\pm\sqrt2)^j
\\&+(2\pm\sqrt2)\(\sum_{k=0}^{p-1}\f{\bi{2k}k^2}{32^k}(1\pm\sqrt2)^k\)^2\pmod{p^2}.
\endalign$$
This, together with the van Hammer-Mortenson congruence
$$\sum_{k=0}^{p-1}(4k+1)\f{\bi{2k}k^3}{(-64)^k}\eq p\l(\f{-1}p\r)\pmod{p^3}$$
(which was conjectured by van Hamme [vH] and proved by E. Mortenson, and further refined by the author [Su6])
 and (3.1), yields that
$$\align&\sum_{k=0}^{p-1}\f{k\bi{2k}k^2}{32^k}(1\pm\sqrt2)^k
\\\eq&\f{2\pm\sqrt2}4(\mp1)^ni^y\f{p(-1)^n-(-1)^y\bar\pi^2}{\bar\pi}
\\=&\f{2\pm\sqrt2}4(\pm1)^ni^y\l(\f p{\bar\pi}-\bar\pi\r)=\f{2\pm\sqrt2}2(\pm1)^ni^yy\sqrt{-2}\pmod{p^2}.
\endalign$$
Therefore
$$\align 2\sqrt2\sum_{k=0}^{p-1}\f{k\bi{2k}k^2}{32^k}P_k=&\sum_{k=0}^{p-1}\f{k\bi{2k}k^2}{32^k}\l((1+\sqrt2)^k-(1-\sqrt2)^k\r)
\\\eq&i^yy\sqrt{-2}\l(\f{2+\sqrt2}2-(-1)^n\f{2-\sqrt2}2\r)\pmod{p^2}
\endalign$$
and hence
$$\sum_{k=0}^{p-1}\f{k\bi{2k}k^2}{32^k}P_k\eq\cases i^yy\sqrt{-2}/2\pmod{p^2}&\t{if}\ p\eq1\ (\mo\ 8),
\\i^{y+1}y=(-1)^{(y+1)/2}y\pmod{p^2}&\t{if}\ p\eq3\ (\mo\ 8).
\endcases$$
Note that if $p\eq1\ (\mo\ 8)$ then
$$i^y=(-1)^{y/2}=(-1)^{(p-1)/8+(x-1)/4}$$
and
$$y\sqrt{-2}=\f p{x-y\sqrt{-2}}-x\eq\f p{2x}-x\pmod{p^2}.$$
Observe that
$$\align&\sum_{k=0}^{p-1}\f{k\bi{2k}k^2}{32^k}Q_k
\\\eq&i^yy\sqrt{-2}\l(\f{2+\sqrt2}2+(-1)^n\f{2-\sqrt2}2\r)
\\\eq&\cases i^yy2\sqrt{-2}\eq4\sum_{k=0}^{p-1}kP_k\bi{2k}k^2/32^k\pmod{p^2}&\t{if}\ p\eq1\ (\mo\ 8),
\\i^{y+1}2y=(-1)^{(y+1)/2}2y\pmod{p^2}&\t{if}\ p\eq3\ (\mo\ 8).\endcases
\endalign$$

Combining the above we immediately get all the results stated in Theorem 1.1. \qed

\medskip
\noindent{\it Proof of Theorem 1.2}. Let $\sqrt{-3}$ be a root of $z^2+3=0$ in the ring $\Z_p$ of $p$-adic integers.
Without loss of generality we suppose that $\pi=x+\sqrt{-3}y\eq0\ (\mo\ p)$ and $x+y\eq1\ (\mo\ 4)$. In view of Lemmas 2.2 and 2.4, we have
$$\sum_{k=0}^{p-1}\f{\bi{2k}k^2}{(-16)^k}\l(\f{-1\pm\sqrt{-3}}2\r)^k\eq (\pm1)^nP_n(\sqrt{-3})
\eq (\pm1)^n(-1)^y\bar \pi \pmod{p^2}.$$ By Lemma 2.3,
$$\bar\pi \eq 2x-\f p{2x}\eq-\f{\sqrt{-3}}2\l(4y-\f p{3y}\r)\pmod{p^2}.$$
Now, applying Lemma 2.1 with $m=16$, $m_*=8(1\pm\sqrt{-3})$, $a=6$ and $b=3\pm\sqrt{-3}$, and noting that
$$\sum_{k=0}^{p-1}(3k+1)\f{\bi{2k}k^3}{16^k}\eq p\ (\mo\ p^2)$$
by [GZ], we get
$$\align\sum_{k=0}^{p-1}\f{k\bi{2k}k^2}{(-16)^k}\l(\f{-1\pm\sqrt{-3}}2\r)^k\eq&\f{3\pm\sqrt{-3}}{12}(\pm1)^n(-1)^y\f{p-\bar\pi^2}{\bar\pi}
\\=&\f{3\pm\sqrt{-3}}6(\pm1)^n(-1)^yy\sqrt{-3}\pmod{p^2}.\endalign$$

Note that $\omega=(-1+\sqrt{-3})/2$ and $\bar\omega=(-1-\sqrt{-3})/2$ are primitive cubic roots of unity.
As $\omega+\bar\omega=-1$ and $\omega\bar\omega=1$, we see that
$$u_n(-1,1)=\f{\omega^n-\bar\omega^n}{\sqrt{-3}}=\l(\f n3\r)\ \t{and}\ v_n(-1,1)=\omega^n+\bar\omega^n=3[3\mid n]-1.$$

Combining the above, similar to the proof of Theorem 1.1, we get all the results in Theorem 1.2. \qed

\medskip
\noindent{\it Proof of Theorem 1.3}. Let $\sqrt{-7}$ be a root of $z^2+7=0$ in the ring $\Z_p$ of $p$-adic integers.
Without loss of generality we suppose that $\pi=x+\sqrt{-7}y\eq0\ (\mo\ p)$ and $x+y\eq1\ (\mo\ 4)$. In view of Lemmas 2.2 and 2.4, we have
$$\sum_{k=0}^{p-1}\f{\bi{2k}k^2}{16^k}\l(\f{1\pm3\sqrt{-7}}2\r)^k\eq (\mp1)^nP_n(3\sqrt{-7})
\eq (\pm1)^n(-1)^y\bar\pi \pmod{p^2}.$$ By Lemma 2.3,
$$\bar\pi \eq 2x-\f p{2x}\eq-\f{\sqrt{-7}}2\l(4y-\f p{7y}\r)\pmod{p^2}.$$
Now, applying Lemma 2.1 with $m=1$, $m_*=(1\pm3\sqrt{-7})/2$, $a=28$ and $b=21\pm\sqrt{-7}$, and noting that
$$\sum_{k=0}^{p-1}(21k+8)\bi{2k}k^3\eq 8p\ (\mo\ p^3)$$
by [Su2], we get
$$\align\sum_{k=0}^{p-1}\f{k\bi{2k}k^2}{16^k}\l(\f{1\pm3\sqrt{-7}}2\r)^k\eq&(\pm1)^n(-1)^y\f{21\mp\sqrt{-7}}{56}\cdot\f{p-\bar\pi^2}{\bar\pi}
\\=&\f{21\mp\sqrt{-7}}{28}(\pm1)^n(-1)^yy\sqrt{-7}\pmod{p^2}.\endalign$$
Note that
$$y\sqrt{-7}=\f p{x-y\sqrt{-7}}-x\eq\f p{2x}-x\pmod{p^2}.$$
Also, $(-1)^y=(\f2p)$ if $p\eq1\ (\mo\ 4)$, and $(-1)^{(y-1)/2}=(\f2p)$ if $p\eq3\ (\mo\ 4)$.
As
$$3\sqrt{-7}u_k(1,16)=\l(\f{1+3\sqrt{-7}}2\r)^k-\l(\f{1-3\sqrt{-7}}2\r)^k$$
and $$v_k(1,16)=\l(\f{1+3\sqrt{-7}}2\r)^k+\l(\f{1-3\sqrt{-7}}2\r)^k,$$
we immediately obtain the desired results from the above.
\qed

\medskip
\noindent{\it Proof of Theorem 1.4}. The proof is similar to the proofs of Theorems 1.1 and 1.2. We apply Lemma 2.1 with $m=256$,
$m_*=64/(2\pm\sqrt3)$, $a=3$ and $b=3\pm2\sqrt3$, and use the congruence
$$\sum_{k=0}^{(p-1)/2}\f{6k+1}{256^k}\bi{2k}k^3\eq p\l(\f{-1}p\r)\pmod{p^4}$$
conjectured by van Hamme [vH] and confirmed by L. Long [L]. \qed

\heading{3. Proof of Theorem 1.5}\endheading

\proclaim{Lemma 4.1} Let $p=2n+1$ be an odd prime. For each $k=0,\ldots,n$ we have
$$\bi{2n-k}k\eq (-1)^k\bi{2k}k(1-p(H_{2k}-H_k))\pmod{p^2}.\tag4.1$$
\endproclaim
\Proof. Observe that
$$\align (-1)^k\bi{2n-k}k=&(-1)^k\prod_{0<j\ls k}\f{p-k-j}j=\prod_{0<j\ls k}\f{j+k-p}j
\\=&\prod_{0<j\ls k}\f{j+k}j\times \prod_{0<j\ls k}\l(1-\f p{j+k}\r)
\\\eq&\bi{2k}k\(1-\sum_{0<j\ls k}\f p{j+k}\)=\bi{2k}k(1-p(H_{2k}-H_k))\pmod{p^2}.
\endalign$$
This concludes the proof. \qed

\proclaim{Theorem 4.1} Let $p=2n+1$ be an odd prime and let $m$ be
an integer not divisible by $p$. Let $h\in\Z^+$ and
$P(x)\in\Z[x]$. Then
$$\aligned&\l(\f{(-1)^hm}p\r)\sum_{k=0}^nP(k)\f{\bi{2k}k^h}{m^k}
\\\eq&\sum_{k=0}^n\f{\bi{2k}k^h}{\bar m^k}\l(\l(1+\f p2q_p(\bar m)\r)P\l(-k-\f12\r)+\f p2P'\l(-k-\f12\r)\r)
\\&-ph\sum_{k=0}^n\f{\bi{2k}k^h}{\bar m^k}\l(P\l(-k-\f12\r)(H_{2k}-H_k)\r)\pmod{p^2},
\endaligned$$
where $\bar m=16^h/m$ and $q_p(\bar m)=(\bar m^{p-1}-1)/p$.
\endproclaim
\Proof. Since
$$\bi{n+k}{n-k}=\bi{n+k}{2k}\eq\f{\bi{2k}k}{(-16)^k}\pmod{p^2}\quad\t{for all}\
k=0,\ldots,n,$$ we have
$$\align\sum_{k=0}^nP(k)\f{\bi{2k}k^h}{m^k}
\eq&\sum_{k=0}^nP(k)\bi{n+k}{n-k}^h\f{(-16)^{hk}}{m^k}
\\=&\sum_{k=0}^nP(n-k)\bi{2n-k}k^h((-1)^h\bar m)^{n-k}\pmod{p^2}.
\endalign$$
So, with the help of Lemma 4.1, we get
$$\sum_{k=0}^nP(k)\f{\bi{2k}k^h}{m^k}\eq\sum_{k=0}^nP(n-k)\f{\bi{2k}k^h}{\bar m^k}(1-p(H_{2k}-H_k))^h((-1)^h\bar m)^n\pmod{p^2}.\tag4.2$$
Since
$$\l(\f p2+x\r)^j-x^j\eq j\f p2x^{j-1}\eq\f p2\f{d}{dx}x^j\ (\mo\
p^2)$$ for all $j=0,1,2,\ldots$, we see that
$$P\l(\f p2+x\r)-P(x)\eq\f p2P'(x)\pmod{p^2}.$$
Thus
$$P(n-k)=P\l(\f p2-k-\f12\r)\eq P\l(-k-\f12\r)+\f p2P'\l(-k-\f12\r)\pmod{p^2}$$
for any $k=0,1,\ldots,n$. So, from (4.2) we obtain
$$\aligned&\l(\f{(-1)^h}p\r)\sum_{k=0}^nP(k)\f{\bi{2k}k^h}{m^k}
\\\eq&\sum_{k=0}^n\l(P\l(-k-\f12\r)+\f p2P'\l(-k-\f12\r)\r)
\f{\bi{2k}k^h}{\bar m^k}(1-ph(H_{2k}-H_k))\bar m^n\pmod{p^2}.
\endaligned\tag4.3$$
Since
$$\bar m^n=\f{4^{(p-1)h}}{m^n}\eq\l(\f mp\r)\pmod p,$$
we have
$$\align \bar m^{p-1}-1=&\l(\bar m^n+\l(\f mp\r)\r)\l(\bar m^n-\l(\f mp\r)\r)
\\\eq&2\l(\f mp\r)\l(\bar m^n-\l(\f mp\r)\r)
\pmod{p^2}
\endalign$$
and hence
$$1+\f p2q_p(\bar m)\eq\l(\f{m}p\r)\bar m^n\pmod{p^2}.$$
Combining this with (4.3) we immediately obtain the desired result. \qed

\medskip
\noindent{\it Proof of Theorem} 1.5. Let $n=(p-1)/2$. Note that $p\mid\bi{2k}k$ for all $k=n+1,\ldots,p-1$.
Applying Theorem 4.1 with  $h=3$ and $P(x)=ax+b$ we then get the desired result.

\medskip

Below we give some corollaries of Theorem 1.5.

\proclaim{Corollary 4.1} Let $p>3$ be a prime. Then
$$\aligned&\l(\f{-1}p\r)\sum_{k=0}^{p-1}\f{k\bi{2k}k^3}{(-64)^k}(H_{2k}-H_k)-\f16
\\\eq&\cases-(3q_p(2)+2)x^2/3\ (\mo\ p)&\t{if}\ p=x^2+2y^2,
\\0\ (\mo\ p)&\t{if}\ p\eq5,7\ (\mo\ 8).\endcases
\endaligned\tag4.4$$
When $p\eq1\ (\mo\ 4)$ and $p=x^2+y^2$ with $x$ odd and $y$ even, we have
$$\sum_{k=0}^{p-1}\f{k\bi{2k}k^3}{64^k}(H_{2k}-H_k)\eq-\f{x^2}3(3q_p(2)+2)\pmod{p}.\tag4.5$$
\endproclaim
\Proof. Recall the known congruence
$$\sum_{k=0}^{p-1}\f{4k+1}{(-64)^k}\bi{2k}k^3\eq p\l(\f{-1}p\r)\pmod {p^3}$$
due van Hamme [vH] and Mortenson.
In view of this and (1.3) and (1.9), Theorem 1.1 with $m=-64$, $a=4$ and $b=1$ yields (4.4).

Let $p\eq1\ (\mo\ 4)$ and write $p=x^2+y^2$ with $x$ odd and $y$ even. The author [Su2, Conj. 5.9] conjectured that
$$\sum_{k=0}^{p-1}\f{4k+1}{64^k}\bi{2k}k^3\eq0\pmod{p^2},$$
and this was confirmed by Z.-H. Sun [S2]. In view of this and (1.0) and (1.9), we obtain (4.5) from Theorem 1.1 with $m=64$, $a=4$ and $b=1$.

\proclaim{Corollary 4.2} Let $p>3$ be a prime. Then
$$\aligned&\sum_{k=0}^{p-1}\f{k\bi{2k}k^3}{16^k}(H_{2k}-H_k)-\f16
\\\eq&\l(\f{-1}p\r)\sum_{k=0}^{p-1}\f{k\bi{2k}k^3}{256^k}(H_{2k}-H_k)-\f16
\\\eq&\cases-2x^2(4q_p(2)+3)/9\ (\mo\ p)&\t{if}\ p=x^2+3y^2,
\\0\ (\mo\ p)&\t{if}\ p\eq2\ (\mo\ 3).\endcases\endaligned\tag4.6$$
\endproclaim
\Proof. van Hamme [vH] conjectured the congruence
$$\sum_{k=0}^{(p-1)/2}\f{6k+1}{256^k}\bi{2k}k^3\eq p\l(\f{-1}p\r)\pmod{p^4}$$
and this was confirmed by L. Long [L]. The author [Su2] conjectured
$$\sum_{k=0}^{(p-1)/2}\f{3k+1}{16^k}\bi{2k}k^3\eq p+2\l(\f{-1}p\r)p^3E_{p-3}\pmod{p^4}$$
(where $E_0,E_1,\ldots$ are Euler numbers), and its mod $p^3$ was proved in [GZ].
In view of this and (1.4)-(1.5) and (1.9), we deduce (4.6) by applying Theorem 1.1 with $(a,b,m)=(6,1,256),(3,1,16)$. \qed

\proclaim{Corollary 4.3} Let $p>3$ be a prime. Then
$$\align&\sum_{k=0}^{p-1}k\bi{2k}k^3(H_{2k}-H_k)-\f16
\\\eq&\cases-2x^2/3\ (\mo\ p)&\t{if}\ (\f p7)=1\ \&\ p=x^2+7y^2,
\\0\ (\mo\ p)&\t{if}\ (\f p7)=-1.\endcases
\endalign$$
Also, we have
$$\aligned&\l(\f{-1}p\r)\sum_{k=0}^{p-1}\f{k\bi{2k}k^3}{4096^k}(H_{2k}-H_k)-\f16
\\\eq&\cases-2(10q_p(2)+7)x^2/21\ (\mo\ p)&\t{if}\ (\f p7)=1\ \&\ p=x^2+7y^2,
\\0\ (\mo\ p)&\t{if}\ (\f p7)=-1.\endcases
\endaligned$$
\endproclaim
\Proof. By [Su2, Theorem 3], we have
$$\sum_{k=0}^{p-1}(21k+8)\bi{2k}k^3\eq8p\pmod{p^4}.$$
Also, van Hamme's conjectural congruence
$$\sum_{k=0}^{(p-1)/2}\f{42k+5}{4096^k}\bi{2k}k^3\eq 5p\l(\f{-1}p\r)\pmod{p^4}$$
was recently confirmed in [OZ].
So, by applying Theorem 1.1 with $(a,b,m)=(21,8,1),(42,5,4096)$ and noting (1.1),(1.6) and (1.9)
we obtain the desired result. \qed

\proclaim{Corollary 4.4} Let $p>3$ be a prime. Then
$$\aligned&\sum_{k=0}^{p-1}\f{k\bi{2k}k^3}{(-8)^k}(H_{2k}-H_k)-\f{(-1)^{(p-1)/2}}6
\\\eq&\cases-2(q_p(2)+1)x^2/3\ (\mo\ p)&\t{if}\ p=x^2+y^2\ (2\nmid x),
\\0\ (\mo\ p)&\t{if}\ p\eq3\ (\mo\ 4).\endcases
\endaligned$$
We also have
$$\aligned&\l(\f{-2}p\r)\sum_{k=0}^{p-1}\f{k\bi{2k}k^3}{(-512)^k}(H_{2k}-H_k)-\f16
\\\eq&\cases-(3q_p(2)+2)x^2/3\ (\mo\ p)&\t{if}\ p=x^2+y^2\ (2\nmid x),
\\0\ (\mo\ p)&\t{if}\ p\eq3\ (\mo\ 4).\endcases
\endaligned$$
\endproclaim
\Proof. The author [Su2] conjectured that
$$\sum_{k=0}^{p-1}\f{3k+1}{(-8)^k}\eq p\l(\f{-1}p\r)+p^3E_{p-3}\pmod{p^4},$$
which was recently confirmed in [CXB]. Also, van Hammer's conjectural congruence
$$\sum_{k=0}^{(p-1)/2}\f{6k+1}{(-512)^k}\bi{2k}k^3\eq p\l(\f{-2}p\r)\pmod{p^3}$$
was recently confirmed by H. Swisher [Sw, Theorem 1.3].
So, by applying Theorem 1.4 with $(a,b,m)=(6,1,-512),(3,1,-8)$ and noting (1.2) and (1.9)
we obtain the desired result. \qed

Our following conjecture seems difficult.

\proclaim{Conjecture 4.1} Let $p$ be an odd prime.

{\rm (i)} If $p\eq1\ (\mo\ 4)$ then
$$\align\sum_{k=0}^{(p-1)/2}\f{\bi{2k}k^3}{(-8)^k}\sum_{k<j\ls 2k}\f1j
\eq&\f12\sum_{k=0}^{(p-1)/2}\f{\bi{2k}k^3}{64^k}\sum_{k<j\ls 2k}\f1j
\\\eq&\f 13\l(\f 2p\r)\sum_{k=0}^{(p-1)/2}\f{\bi{2k}k^3}{(-512)^k}\sum_{k<j\ls 2k}\f1j\ (\mo\ p^2);
\endalign$$
when $p\eq3\ (\mo\ 4)$ we have
$$\align\sum_{k=0}^{(p-1)/2}\f{\bi{2k}k^3}{(-8)^k}\sum_{k<j\ls 2k}\f1j
\eq&-\f 72\sum_{k=0}^{(p-1)/2}\f{\bi{2k}k^3}{64^k}\sum_{k<j\ls 2k}\f1j\ (\mo\ p^2),
\\\sum_{k=0}^{(p-1)/2}\f{\bi{2k}k^3}{64^k}\sum_{k<j\ls 2k}\f1j
\eq&-\l(\f 2p\r)\sum_{k=0}^{(p-1)/2}\f{\bi{2k}k^3}{(-512)^k}\sum_{k<j\ls 2k}\f1j\ (\mo\ p^2),
\endalign$$

{\rm (ii)} If $p\eq1\ (\mo\ 3)$ then
$$\sum_{k=0}^{(p-1)/2}\f{\bi{2k}k^3}{16^k}\sum_{k<j\ls 2k}\f1j
\eq\f12\l(\f{-1}p\r)\sum_{k=0}^{(p-1)/2}\f{\bi{2k}k^3}{256^k}\sum_{k<j\ls 2k}\f1j\ (\mo\ p^2).$$
If $p\eq 2\ (\mo\ 3)$ then
$$\sum_{k=0}^{(p-1)/2}\f{\bi{2k}k^3}{256^k}\sum_{k<j\ls 2k}\f1j\eq0\ (\mo\ p^2).$$

{\rm (iii)} If $p>3$ and $p\eq3,5,6\pmod7$, then
$$\sum_{k=0}^{(p-1)/2}\bi{2k}k^3\sum_{k<j\ls 2k}\f1j\eq0\pmod{p^2}.$$
\endproclaim

\enddocument